\documentclass{article}
\usepackage{ijcai19}
\usepackage{amsmath}
\usepackage{amssymb}
\usepackage{titlesec}

\usepackage{times}
\usepackage{soul}
\usepackage{url}
\usepackage{graphicx}
\usepackage{booktabs}
\usepackage{color}
\usepackage{algorithm}
\usepackage{algpseudocode}
\input{mysymbol.sty}
\usepackage{theorem}

\newtheorem{thm}{Theorem}
\newtheorem{lemma}{Lemma}

\newtheorem{assumption}{Assumption}

\def \QED {$\blacksquare$}
\def \w {{\mathbf{w}}}
\def \x {{\mathbf{x}}}
\def \z {{\mathbf{z}}}
\def \v {{\mathbf{v}}}
\def \u {{\mathbf{u}}}

\def \e {{\mathbf{e}}}

\providecommand{\norm}[1]{\left\|#1\right\|}

\title{\bf Escaping Saddle Points with the Successive Convex Approximation Algorithm}

\author{
	Amrit Singh Bedi$^1$\footnote{Contact Author}\and
	Ketan Rajawat$^2$\And
	Vaneet Aggarwal$^{3}$\\
\affiliations
$^{1,2}$Dept. of EE, IIT Kanpur, India\\
$^3$Dept. of IE, Purdue University, USA\\
\emails
\{amritbd, ketan\}@iitk.ac.in,
vaneet@purdue.edu
}

\begin{document}
	\maketitle
\begin{abstract}
Optimizing non-convex functions is of primary importance in the vast majority of machine learning algorithms. Even though many gradient descent based algorithms have been studied, successive convex approximation based algorithms have been recently empirically shown to converge faster. However, such successive convex approximation based algorithms can get stuck in a first-order stationary point. To avoid that, we propose an algorithm that perturbs the optimization variable slightly at the appropriate iteration. In addition to achieving the same convergence rate results as the non-perturbed version, we show that the proposed algorithm converges to a second order stationary point. Thus, the proposed algorithm escapes the saddle point efficiently and does not get stuck at the first order saddle points.


\end{abstract}

\section{Introduction}\label{sec:intro}

Optimization, as a field, has found strong applications in a lot of science and technology areas. Some of them include mechanics, economics and management, electromagnetics,  communications, scheduling, control systems, smart grids, and artificial intelligence \cite{jain2017non}.  Most of the algorithms in supervised and unsupervised machine learning involve solving optimization problems. Most of these problems are non-convex optimization problems such as matrix completion \cite{Matrix}, phase retrieval \cite{Sun2018}, tensor decomposition \cite{Ge15},  tensor completion \cite{liu2016low,wang2017efficient}, dimensionality reduction \cite{wang2018tensor}, etc., and efficient mechanisms to solve the problems is an important research direction.  Even though many of these problems are non-convex, a popular
workaround is to relax the problems to a convex optimization problem.  However, this approach may be lossy and nevertheless presents significant challenges for large scale optimization. Thus, novel algorithms for solving non-convex optimization problems are of interest in machine learning, and this paper focuses on proposing and analyzing a novel algorithm for solving non-convex optimization problems.

Efficient approaches for solving non-convex optimization problems have been widely studied. Most common approaches use a version of gradient (or sub-gradient) descent \cite{jin,7029630,7536166}. In contrast, the authors of \cite{6675875,scutari2017parallel,liu2018stochastic} empirically demonstrated that successive convex approximation (SCA) based approaches can outperform the gradient descent based approaches, even though both obtain the same convergence rate. Thus, we consider an algorithm based on the SCA approach. Recently, \cite{jin} have proposed an approach in which a perturbed form of gradient descent algorithm is shown to converge to a second-order stationary point, and can escape the saddle points. However, there has been no such study for the SCA approach for convergence to second order stationary point (to the best of our knowledge), even though the SCA approaches are more efficient. This paper tackles this problem and shows that the  SCA can also avoid the saddle points. Further, the convergence rate guarantees are maintained with perturbation. 

To be specific, in this paper we are interested in the question that whether SCA algorithm converges to a local minima in the number of iterations which are almost dimension free? The same question is answered for the gradient descent algorithm in \cite{jin}. We extend the analysis to SCA algorithms and show that a perturbed successive convex approximation (P-SCA) algorithm converges to $\epsilon$ second order stationary point in $\mathcal{O}\left(\frac{1}{\epsilon^2}\right)$ with an extra multiplicative factor of \textit{polylog($d$)} where $d$ is the underlying dimension. Note that the dependence on the underlying dimension $d$ is better than the $\Omega(d^4)$ dependence obtained in \cite{Ge15}. The convergence rate obtained in this paper is similar to the well-known rate of gradient decent algorithm convergence to a first order stationary point (FOSP) with an additional log factor \cite{nesterov1998introductory}. This comes from the structure of the proposed SCA type algorithm. It is emphasized that the proposed algorithm does not exhibit a gradient descent like update, and hence the analysis of \cite{jin} cannot be directly extended. To do so, we show that the proposed algorithm is an inexact version of the standard gradient descent algorithm, which is then used to derive the convergence rate result \cite{bertsekas2000gradient,so2017non}. We note that the results in this paper can be generalized to an arbitrary inexact gradient descent algorithm with an additional assumption that the norm of the error in the gradient is bounded.

\section{Related Work}

In literature, gradient/subgradient descent based approaches are popular for solving the optimization problems in an iterative manner \cite{spall2005introduction}. To handle the non-convexity in the objective function, gradient averaging based methods are proposed in \cite{Ruszczyński1980,JMLR:v15:bach14a}. The other related works in \cite{netrapalli2013phase,7029630,Matrix} require smart initialization and have demonstrated faster local convergence. 

Apart form the first order methods mentioned above, there are methods such as Majorization-Minimization (MM) \cite{7547360,mairal2013stochastic} and SCA \cite{6675875} which are used for the non-convex optimization. Both of these methods utilize the same technique of solving a sequence of approximate surrogate functions under different conditions. At the point of approximation, MM requires the surrogate function to upper bound the non-convex objective function, while for SCA the surrogate function needs to be convex and the upper bound condition is not required. Therefore, SCA helps to better approximate the non-convex objective and results in better convergence speeds \cite{liu2018stochastic}. 

For the approximation based approaches to optimize non-convex objectives, the convergence to FOSP is discussed in literature \cite{scutari2017parallel,liu2018stochastic}. In this paper, we are interested in making these approximation based technique to escape saddle points and converge to the second order stationary point (SOSP). 

To escape the saddle point, two popular techniques namely perturbation based \cite{jin} and second-order information based methods \cite{Nesterov2006} have been used. A cubic regularization based algorithm is proposed in  \cite{Nesterov2006} with convergence to $\epsilon$ SOSP in $\mathcal{O}(\frac{1}{\epsilon^{3/2}})$ number of iterations. In practice, the second order information based methods require Hessian inversion at each algorithm iteration, which is computationally expensive. On the other hand, in perturbation based methods, noise from a given distribution (uniform in this paper) is added to the algorithm  update (when in the neighborhood of a saddle point) to escape the saddle point which is computationally inexpensive. This work aims to achieve advantages of both the SCA based methods equipped with perturbation to escape the saddle points efficiently and converge to $\epsilon$ SOSP.

	\section{Preliminaries}
	
	In this paper, we are solving a non-convex optimization problem with the help of gradient based methods. To follow the analysis, we describe some useful definitions which will be referred  throughout the paper. 
	
	\textbf{Definition 1} For a differentiable function $U(\x)$, a point $\x$ is $\epsilon$ first order stationary point (FOSP) if $\norm{\nabla U(\x)}\leq \epsilon$. 
	
	\textbf{Definition 2} (Strict saddle point) For a twice differentiable function, a point $\x$ is a strict saddle point if $\x$ is a FOSP and $\lambda_{\min}(\nabla^2 U(\x))<0$.  For a strict saddle point, the second order stationary point (SOCP) is  a local minima of the non-convex objective as defined next. 
	
	\textbf{Definition 3} ($\epsilon$  second order  order stationarity)  For a non-convex function $U(\x)$, a point $\x$ is $\epsilon$ second order stationary  if 
	\begin{align}\label{bounds}
	\norm{\nabla U(\x)} \leq \epsilon\ \ \ \ \text{and} \ \ \ \ \lambda_{\min}(\nabla^2 U(\x))\geq -\sqrt{L_{2}\epsilon}. 
	\end{align}where $L_2$ is Hessian Lipschitz constant. Note that some other parameters can also be used to define the upper and lower bounds in \eqref{bounds} but the mentioned values are used to make the first condition for SOSP to be same as that for FOSP as in \cite{jin}.  
	
	\section{Problem Formulation}\label{sec:prob0}
	In this work, we focus on the minimization of a non-convex objective function $U(\x)$
	given by
	\begin{align}\label{eq:main_prob_uncons}
	\min_{\bbx}\ \ & U(\bbx)
	\end{align}
	where $\x\in\mathbb{R}^n$ is the optimization variable, $U(\x):\mathbb{R}^n\rightarrow\mathbb{R}$ is a twice differentiable non-convex function. For a general non-convex optimization problem, it is well known that an algorithm cannot achieve global minima using first order algorithms \cite{Ge15,jin}. In literature, there are various algorithms proposed which converges to the approximate FOSP of the problem \eqref{eq:main_prob_uncons}. A FOSP could be a local minima, local maxima, or a saddle point. It is sufficient to achieve convergence to FOSP in convex settings because it coincides with the global minima of the function. However, for the non-convex functions, it may be highly suboptimal to achieve FOSP if it denotes the local maxima and we are minimizing the function and vice versa. In addition, the converged FOSP can be a saddle point for the unconstrained optimization problem and algorithm might get stuck there. 
	Hence, it becomes important to look for the convergence to a second order stationary point. In literature, various algorithms are proposed for problem \eqref{eq:main_prob_uncons} with the guarantee of convergence to a second order stationary point \cite{curtis2017trust,jin,levy2016power}. However, they are based on the standard gradient descent algorithm. 
	
	Recently,  an SCA based algorithm is proposed by ~\cite{scutari2017parallel}  to solve the problem of \eqref{eq:main_prob_uncons} in an iterative manner. We utilize the similar ideas and use SCA to solve the unconstrained version of the problem.  The SCA algorithm is shown to outperform the other first order methods in  literature \cite{scutari2017parallel}. In this method, the non-convex objective function $U(\x)$ is approximated by a surrogate convex function $\tilde{U}(\x;\mathbf{y})$ approximated at $\mathbf{y}$. Utilizing this idea, the standard SCA algorithm to solve the problem in \eqref{eq:main_prob_uncons} is summarized in Algorithm \ref{algo1}. 
	
	\begin{algorithm}
		\caption{:Successive convex approximation (SCA) }\label{algo1}
		\begin{algorithmic}[1]
			\State {\bf Set} $t=0$, initialize $\eta=(0,1]$,   
			\State \textbf{STOP} if $\bbx_t$ is a FOSP
			\State Solve the following optimization problem to get $\hat{\bbx}(\bbx_t)$
			\begin{align}\label{eq:main_prob_approx_algo}
			\hat{\bbx}(\bbx_t):=\ \ \arg\min_{\bbx\in\mathbb{R}^n}\ \ & \tilde{U}(\bbx;\bbx_{t}) ,
			\end{align}
			\State Set $\bbx_{t+1}=\bbx_t+\eta(\hat{\bbx}(\bbx_t)-\bbx_t)$
			\State Set $t=t+1$ and go to step 2
		\end{algorithmic}
	\end{algorithm}
	
	In Algorithm \ref{algo1}, at each step $t$, a convex approximation $\tilde{U}(\x;\x_t)$ of the non-convex objective function $U(\x)$ is used to perform the update in step 3. This algorithm converges to the first order  stationary point as derived in \cite[Theorem 2]{scutari2017parallel}, which states that for a constant step size $\eta$, the sequence $\{\x_t\}$ is bounded and each of its limit points is stationary point of the optimization problem in \eqref{eq:main_prob_uncons}. But the results presented in  \cite{scutari2017parallel} are asymptotic in nature and does not provide the number of iterations required to converge to an $\epsilon$ optimal solution. This results is important because it characterizes the convergence behavior of the proposed algorithm.  Moreover, as stated earlier, the convergence to FOSP is not a sufficient condition of convergence to local minima. This is because a FOSP could be a saddle point and the proposed algorithm might get stuck there. 
	
	In this paper, we are interested in deriving the number of iterations required to achieve an $\epsilon$ FOSP. In addition to that, we want to characterize the number of iterations required to achieve an $\epsilon$ second order stationary point. 	In literature there are different techniques proposed to escape the saddle points and converge to SOSP. For instant, some of the approached include methods utilizing second order information \cite{nesterov1998introductory}, and random perturbation based methods \cite{jin}. Motivated from the advantages of perturbation based methods as discussed in \cite{jin} and advantages of SCA methods \cite{scutari2017parallel}, we use the perturbed method to make the SCA algorithm escape saddle points effectively. A perturbed successive convex approximation (P-SCA) algorithm is summarized in Algorithm \ref{algo3}.
	
	\begin{algorithm}
		\caption{:P-SCA Algorithm }\label{algo3}
		\begin{algorithmic}[1]
			\State {\bf Set} $t=0$, initialize $\chi\leftarrow3\max\{\log\left(
			\frac{dL_1\triangle_U}{c\epsilon^2\delta},4\right)\}$, $\eta\leftarrow\frac{c}{L_1}$, $r\leftarrow\frac{\epsilon\sqrt{c}}{L_1\chi^2}$, $g_{\text{th}}\leftarrow\frac{\epsilon\sqrt{c}}{\chi^2}$, $f_{\text{th}}\leftarrow\frac{c}{\chi^3}\sqrt{\frac{\epsilon^3}{L_2}}$, $t_{\text{th}}=\frac{(1-s)\chi L_1}{c^2\sqrt{\epsilon L_2}}$, $t_{\text{noise}}=-t_{\text{th}}-1$, $s\in(0,1)$, and $\bbx_0$  
			
			\State If $\norm{\nabla U{(\x_t)}}\leq g_{\text{th}}$ and $t-t_{\text{noise}}>t_{\text{th}}$ then 
			
			$\tilde{\x}_t\leftarrow \x_t$, and $t_{\text{noise}}\leftarrow t$
			
			$\x_t\leftarrow\tilde{\x}_t+\bbxi_t$, where $\bbxi$ is uniformly $\sim \mathbb{B}_0(r)$ with radius $r$ 
			
			\State If $t-t_{\text{noise}}=t_{\text{th}}$ and $U(\x_t)-U(\tilde{\x}_{t_{\text{noise}}})>-(1-s)f_\text{th}$ then 
			
			\textbf{return} $\tilde{\x}_{t_{\text{noise}}}$
			
			\State Solve the following optimization problem to get $\hat{\bbx}(\bbx_t)$
			\begin{align}\label{eq:main_prob_approx_algo_no_constraints}
			\hat{\bbx}(\bbx_t):=\ \ \arg\min_{\bbx}\ \ & \tilde{U}(\bbx;\bbx_{t}) 
			\end{align}
			\State Set $\bbx_{t+1}=\bbx_t+\eta(\hat{\bbx}(\bbx_t)-\bbx_t)$
			\State Set $t=t+1$ and go to Step 2
		\end{algorithmic}
	\end{algorithm}
	
	In Algorithm \ref{algo3}, the standard SCA algorithm updates are used to reach FOSP. If the condition for FOSP are satisfied (which is nothing but the gradient norm is smaller than a threshold value $g_\text{th}$), a random perturbation from uniform ball $\mathbb{B}_0(r)$ of radius $r$ around $\tilde{\x}_t$ is added. This happens at most once every $t_\text{th}$ number of iterations. After adding perturbation, if the function value is not decreased by the threshold $(1-s)f_\text{th}$ (where $s\in(0,1)$ is a constant and defined later in the paper) then $\tilde{\x}_{t_{\text{noise}}}$ is returned as the output of Algorithm \ref{algo3}. Here $t_{\text{noise}}$ denotes the last iteration at which the optimization variable was perturbed with respect to current iteration of the algorithm. It is proved in \cite{scutari2017parallel} that after the addition of this random perturbations, the standard gradient descent iterates escapes the saddle points efficiently. We extend the analysis in \cite{scutari2017parallel} and \cite{jin} papers to prove that the proposed P-SCA algorithm also escapes the saddle points in a similar manner and $\tilde{\x}_{\text{noise}}$ will be in the neighborhood of $\epsilon$ SOSP. We note that the steps in Algorithm \ref{algo3} do not exhibit a gradient descent like update. Hence, it is difficult to analyze the dynamics of the proposed algorithm. To perform the analysis, it becomes important to look at the proposed algorithm from the perspective of inexact gradient descent algorithms \cite{inexact,8268092,8598992,so2017non}. The meaning of term `inexact' means that the gradient value used for the algorithm update is in error. The error may arise due to stochastic nature of the gradient or partial availability of the gradient \cite{bertsekas2000gradient} etc. In the above mentioned references, all the inexact gradient algorithms are considered for convex objective function settings. The analysis for asymptotic convergence to FOSP under non-convex objective is performed in \cite{bertsekas2000gradient} but with diminishing step size. The work in this paper focuses on the more general scenario of non-convex objective function and inexactness of the gradient \cite{xiao2014proximal}.
	
	Hence, the main idea for the inexact version is to interpret the algorithm as a gradient descent algorithm with error. It is proved in literature that under some regularity conditions on the error term in the gradient, the algorithm converges to (or provide close approximation to) the optimal value under convexity assumption . In order to proceed with analysis for the proposed P-SCA algorithm, we first show that the SCA algorithm is nothing but an inexact version of the standard gradient descent algorithm. From step 5 of the Algorithm \ref{algo3}, we have
	\begin{align}
	\bbx_{t+1}=\bbx_t+\eta(\hat{\bbx}(\bbx_t)-\bbx_t),
	\end{align}
	Next, add and subtract the original gradient $\nabla U(\x_t)$ as follows
	\begin{align}
	\bbx_{t+1}=\bbx_t+\eta(\nabla U(\x_t)-\nabla U(\x_t)+\hat{\bbx}(\bbx_t)-\bbx_t).
	\end{align}
	Finally, we can write the P-SCA algorithm update in the following alternative form
	\begin{align}\label{SCA-GD}
	\bbx_{t+1}=\underbrace{\bbx_t-\eta\nabla U(\x_t)}_{\text{GD  update}}+\eta\underbrace{\left[\nabla {U}(\x_t)+\bbx_t-\hat{\bbx}(\bbx_t)\right]}_{:=\bbe_t}.
	\end{align}
	 The first term in \eqref{SCA-GD} is similar to the standard gradient descent algorithm. The second term in \eqref{SCA-GD} corresponds to the error in the gradient and denoted by $\e_t$. Hence the final gradient used for the update is $\nabla U(\x_t)+\e_t$ as follows
	 \begin{align}\label{SCA-GD2}
	 	\bbx_{t+1}={\bbx_t-\eta(\nabla U(\x_t)+\e_t)}.
	 	\end{align}
	For the further analysis in this paper, we will utilize the algorithm update form as presented in \eqref{SCA-GD}.
	
	\section{Convergence analysis} 
	This sections details about the convergence rate analysis of the proposed P-SCA algorithm and show that the number of iterations required to converge to a second order stationary point are of the order of $\mathcal{O}\left((\text{polylog}(d))\frac{L_1(U(\x_0)-U(\x^\star))}{\epsilon^2}\right)$. This section also establishes the result that the proposed algorithm indeed escapes the saddle point with probability $1-\delta$ for any $\delta>0$. Note that the updates in the proposed algorithm utilize a convex approximations $\tilde{U}(\x;\x_t)$ at each step $t$  for the non-convex function $U(\x)$ \cite{scutari2017parallel}. Before discussing the mathematical results, some assumptions are required to hold for the objective function and its convex approximation at each $t$. All the required assumptions are provided next. That will be followed by the main results.
	
	\subsection{Required assumptions}
	The non-convex objective function is required to satisfy some other technical condition as mentioned here.  
	\begin{assumption}\label{ass_1}
		The continuously differentiable objective function $U(\x)$ and its gradient $\nabla U(\x)$ are Lipschitz continuous with parameter $L_0>0$ and $L_1>0$, respectively. This  implies that
		\begin{align}
		\norm{ U(\x)-U(\bby)}\leq L_0\norm{\x-\bby}
		\end{align}
		and 
		\begin{align}
		\norm{\nabla U(\x)-\nabla U(\bby)}\leq L_1 \norm{\x-\bby}
		\end{align} 
		for all $\x$ and $\bby$. 
	\end{assumption}
	Apart from the objective function, the convex surrogate function needs to satisfy the following assumption. 
	\begin{assumption}\label{ass_2}
		The convex approximation function $\tilde{U}(\x;\bby)$ to non-convex objective function $U(\x)$ at $\bby$ is continuously differentiable with respect to its first argument such that 	
		\begin{itemize}
			\item [B1] $\tilde{U}(\cdot;\bby)$ is uniformly strongly convex with parameter $C$, which implies that 
			\begin{align}
			\!\!\!\!\!\!\!\!\left(\nabla \tilde{U}(\x;\bby)-\nabla \tilde{U}(\bbz;\bby)\right)^T(\bbx-\bbz)\geq C\norm{\bbx-\bbz}^2
			\end{align}
			for all $\x$ and $\bbz$. 
			
			\item[B2] The approximation function gradient $\nabla \tilde{U}(\cdot;\bby)$ is equal to the original function gradient at the approximation point $\bby$ \begin{align}
			\nabla \tilde{U}(\bby;\bby)=\nabla {U}(\bby).
			\end{align}This property is  the key to represent the P-SCA algorithm as an inexact gradient descent algorithm.
		\end{itemize}	
	\end{assumption}
	Another assumption required for the objective function which upper bounds the rate of change of the Hessian.
	\begin{assumption}\label{ass_3}
		The objective function $U(\x)$ is Hessian Lipschitz with parameter $L_{2}$, which implies that
		\begin{align}\label{Hessia_lip}
		\norm{\nabla^2 U(\x)-\nabla^2 U(y)}\leq L_{2} \norm{\x-\bby}
		\end{align}
		for all $\x$ and $\bby$.

	\end{assumption}
	
	All of the above mentioned assumptions are standard and have been considered in literature~\cite{scutari2017parallel}. Assumption \ref{ass_1} states that the non-convex objective function is Lipschitz which means that  the gradient of the objective function is smooth and does not changes arbitrarily. This is an important assumption since it provides an upper bound on the gradient at two difference points in terms of the difference between the points itself. Next, assumption \ref{ass_2} is for the convex approximation $\tilde{U}(\x;\mathbf{y})$ to the non-convex objective function $U(\x)$ at $\mathbf{y}$, and states that $\tilde{U}(\x;\mathbf{y})$ is strongly convex at given $\mathbf{y}$ with parameter $C$. This assumption assures that the convex optimization problem in \eqref{eq:main_prob_approx_algo_no_constraints} has a unique solution. The last assumption \ref{ass_3} states that the gradient of $U(\x)$ is smooth and hence the Hessian of the function does not changes abruptly. All of the above mentioned assumptions are utilized to prove the convergence results presented in the following subsection of this paper. 
	
	\subsection{Main results}
	In this subsection, we show that it is possible for the SCA algorithm to converge to the second order stationary point with a simple modification. In the SCA algorithm updates, once the gradient norm is less than a threshold value $g_{\text{th}}$, we add a small random perturbation to current iterate $\tilde{\x}_t$. In other words, a random perturbation is added to the algorithm updates at most once after every $t_{\text{th}}$ iterations. To make the analysis simple, similar to \cite{jin} we assume that the random perturbation $\bbxi_t$ is uniformly sampled from a ball $\mathbb{B}_0(r)$ of radius $r$ centered around $\tilde{\x}_t$. After adding the perturbation at iteration $t$, if the function value does not decrease by $(1-s)f_\text{th}$ after $t_\text{th}$ iterations, then the algorithm stops and return $\tilde{\x}_{t_{\text{noise}}}$. This paper proves that the output $\tilde{\x}_{t_{\text{noise}}}$ is essentially the second order stationary point for the problem in \eqref{eq:main_prob_uncons}. The main result is presented next in Theorem \ref{thm1}.
	
	\begin{thm}\label{thm1}
		Let the assumptions \ref{ass_1}-\ref{ass_3} are satisfied, there exits a constant $c_{\max}$ such that, for any $\delta>0$, $\epsilon\leq \frac{L_1^2}{L_{2}}$, $\triangle_U\geq U(\x_0)-U(\x^\star)$, and $c\leq c_{\max}$, the output of Algorithm \ref{algo3} is a $\epsilon
		$ second order stationary point with probability {$(1-\delta)$} for problem \eqref{eq:main_prob_uncons}. To achieve the $\epsilon$ second order stationary point, the number of iterations required for the algorithm is given by 
		\begin{align}
		\mathcal{O}\left(\frac{L_1\triangle_U}{\epsilon^2}\log^4\left(\frac{d L_1\triangle_U}{\epsilon^2\delta}\right)\right).
		\end{align}
	\end{thm}
	
	It is interesting to note that the convergence rate of the proposed algorithm remains the same as that of the one in \cite[Theorem 3]{jin}. Theorem \ref{thm1} establishes the fact that the algorithm terminates in finite number of iterations. Before discussing the proof of Theorem \ref{thm1}, there are some other results which we need to discuss first.  The algorithm operation can be divided in two scenarios 1) the iterate is not close to FOSP and the gradient $\nabla U(\x_t)$ is large and 2) when gradient is small but the Hessian $\nabla U(\x_t)$ has a very small minimum eigenvalue. First, to prove the convergence to FOSP, we need to show that a single iteration of the proposed algorithm results in a descent direction for the objective function $U(\x)$ which is stated in Lemma \ref{lem:descent}.

	\begin{lemma}\label{lem:descent}
		Under the Assumptions 1-2, and for a $L_{2}$ Hessian Lipschitz objective function $U(\x)$, one iterate of Algorithm \ref{algo3} is a descent direction, which implies that 
		\begin{align}
		U(\x_{t+1})\leq U(\x_t)-\eta'\norm{(\tilde{\x}(\x_t)-\x_t)}^2
		\end{align} 
		where $\eta':=\eta C-\eta^2\frac{L_{1}}{2}$ with $\eta\leq\frac{2C}{L_{1}}$.
	\end{lemma}
	
	\begin{myproof}
		\normalfont
		Note that $\hat{\x}(\x_t)$ is the solution of strongly convex problem in \eqref{eq:main_prob_approx_algo_no_constraints}, hence from the first order optimality condition, it holds that
		\begin{align}
		(\bby-\hat{\x}(\x_t))^T\nabla \tilde{U}(\hat{\x}(\x_t); \x_t) \geq 0.
		\end{align}
		By selecting $\bby=\x_t$ and add subtract $\nabla \tilde{U}(\x_t; \x_t)$ to the gradient term, we get
		\begin{align}
		(\x_t\!-\!\hat{\x}(\x_t))^T\!\!\left(\!\!\nabla \tilde{U}(\hat{\x}(\x_t); \x_t)\!-\!\nabla \tilde{U}(\x_t; \x_t)\!\!+\!\!\nabla \tilde{U}(\x_t; \x_t)\right)\!\! \geq\!\! 0.
		\end{align}
		Following   the inequality in Assumption, it holds that \ref{ass_2} $\nabla{U}(\x_t)=\nabla\tilde{U}(\x_t; \x_t)$. This implies
		\begin{align}\label{first}
		(\x_t-\hat{\x}(\x_t))^T\nabla {U}(\x_t) \geq C\norm{\hat{\x}(\x_t)-\x_t}^2.
		\end{align}
		From the statement of Assumption \ref{ass_1}, it holds that the objective function  $U(\x)$ is Lipschitz continuous gradient with parameter $L_1$, which implies that 
		\begin{align}
		U(\x_{t+1})\leq & U(\x_t)+\nabla U(\x_t)^T\left(\x_{t+1}-\x_t\right)\nonumber \\&+\frac{L_{1}}{2}\norm{\x_{t+1}-\x_t}^2.
		\end{align} 
		From the update for $\x_{t+1}$ in Algorithm \ref{algo3}, we get
		\begin{align}
		U(\x_{t+1})\leq &U(\x_t)+\eta\nabla U(\x_t)^T\left(\tilde{\x}(\x_t)-\x_t\right)\nonumber\\&+\frac{\eta^2L_{1}}{2}\norm{\tilde{\x}(\x_t)-\x_t}^2.
		\end{align} 
		Applying the upper bound in \eqref{first}, we get
		\begin{align}
		\!\!\!\!U(\x_{t+1})&\leq U(\x_t)-\eta C\norm{\hat{\x}(\x_t)-\x_t}^2\nonumber\\&\hspace{5mm}+\frac{L_{1}}{2}\norm{\eta(\tilde{\x}(\x_t)-\x_t)}^2\\
		&= U(\x_t)+\left(-\eta C+\frac{\eta^2L_{1}}{2}\right)\norm{(\tilde{\x}(\x_t)-\x_t)}^2\nonumber\\
		&= U(\x_t)-\eta'\norm{(\tilde{\x}(\x_t)-\x_t)}^2.\label{descent}
		\end{align} 
		where $\eta':=\eta C-\eta^2\frac{L_{1}}{2}$. Hence, one run of the Algorithm \ref{algo3} results in a function value decrement with $\eta<\frac{2C}{L_{1}}$. 	
	\end{myproof}
	
	The above mentioned analysis helps us to derive an upper bound on the gradient error norm. Consider the additional term in \eqref{SCA-GD}, we have
	\begin{align}
	\norm{\bbe_t}=&\norm{\nabla {U}(\x_t)+\bbx_t-\hat{\bbx}(\bbx_t)}\\
	\leq&\norm{\nabla {U}(\x_t)}+\norm{\hat{\bbx}(\bbx_t)-\bbx_t}.
	\end{align}From the upper bound  in \eqref{first} and Lipschitz property of the objective function, we have
	\begin{align}\label{error_bound}
	\norm{\bbe_t}\leq& L_0\left(1+\frac{1}{C}\right)=:\mathcal{D}.
	\end{align}
	Before proceeding, let us define the following parameters as specified in \cite{scutari2017parallel} which are used for deriving the mathematical results next. 
	\begin{align}\label{eq:definitions}
	& F:=\frac{\eta L_1}{L_{2}^2}\gamma^3\cdot\log^{-3}\left(\frac{d\kappa}{\delta}\right), \  G:=\frac{\sqrt{\eta L_1}}{L_{2}}\gamma^2\cdot \log^{-2}\left(\frac{d\kappa}{\delta}\right)\nonumber\\
	&L:= \sqrt{\eta L_{1}}\frac{\gamma}{L_{2}}\cdot\log^{-1}\left(\frac{d\kappa}{\delta}\right), \ \ \ \ \mathcal{T}:=\frac{\log\left(\frac{d\kappa}{\delta}\right)}{\eta\gamma}
	\end{align}
	where $\gamma$ is the negative eigenvalue and the condition number $\kappa$ is is given by $\kappa=\frac{L_1}{\gamma}\geq 1$. 
	Once the algorithm iterate $\x_t$ is near a FOSP, it holds that  $\norm{\nabla U(\x_t)}\leq g_{\text{th}}$ (gradient is small) and the minimum eigen value of Hessian is largely negative $\lambda_{\min}(\nabla^2U(\x_t))\leq -\sqrt{L_{2}\epsilon}$. When $\x_t$ is near a FOSP,  then we add perturbation to $\x_t$ followed by  SCA updates for $t_{\text{th}}$ steps. We prove that after these $t_\text{th}$ steps, the function value will decrease by at least $f_{\text{th}}$  with high probability. This result is formalized in Lemma \ref{function_value}.
	
	\begin{lemma}\label{function_value}
		Under the Assumption \ref{ass_1}-\ref{ass_2}, for any $\delta\in(0,d\kappa/e]$, if we have $\tilde{\x}$ such that $\norm{\nabla U(\tilde{\x})}\leq G$ (gradient is small) and $\lambda_{\min}(\nabla^2 U(\tilde{\tilde{\x}}))\leq -\gamma$ (sufficiently negative), then if we perform $\x_0=\tilde{\x}+\bbxi$ where $\bbxi\in\mathbb{B}_{\tilde{\x}}(r)$ with radius $r=\frac{L}{\kappa\log\left(\frac{d\kappa}{\delta}\right)}$, it holds that
		\begin{align}
		U(\x_T)-U(\tilde{\x})\leq -F+16L_2\eta^3\mathcal{D}^3
		\end{align}
		for any $T\geq\frac{1}{c_{\max}\mathcal{T}} $with probability $1-\delta$ when the step size is selected as {$\eta\leq\frac{C_{\max}}{L_1}$ and $\delta=\frac{dL_1}{\sqrt{L_2\epsilon}}\exp(-\chi)$}.
	\end{lemma} 
	The result in above lemma states that adding the perturbation decreases the function values further. By selecting $\eta=\frac{c}{L_1}$, $\gamma=\sqrt{L_2\epsilon}$, and $\delta=\frac{dL_1}{\sqrt{L_2\epsilon}}\exp(-\chi)$. We can restate the result in Lemma \ref{function_value} follows. If $\tilde{\x}_t$ is a FOSP, then after adding perturbation $\x_t=\tilde{x}_t+\bbxi_t$, it holds with probability $1-\delta$ that
	\begin{align}
	U(\x_{t+th})-U(\tilde{\x}_t)\leq -f_\text{th}+\frac{16L_2c^3\mathcal{D}^3}{L_1^3}
	\end{align} with $t_\text{th}=\frac{\chi}{c^2}\frac{L_1}{\sqrt{L_2\epsilon}}$.  Next lemma states that addition of perturbation indeed results in escaping the saddle points. 
	
	\begin{lemma}\label{funvtion_bound}
		Under the conditions similar to Lemma \ref{function_value}, let $\bbg_1$ is the minimum eigen vector of $\nabla^2U(\tilde{\x})$. If we consider two sequences $\u_t$ and $\w_t$ such that 
		\begin{align}
		\norm{\u_0-\tilde{\x}}\leq r, \ \ \ \w_0=\u_0+\mu r\bbg_1\ ; \ \mu\in[\delta/(2\sqrt{d}),1],
		\end{align}
		then it holds that
		\begin{align}
		\min&\{U(\u_T)-U(\u_0),U(\w_T)-U(\w_0)\}\nonumber \\ 
		&\leq -2.5F{+16L_2\eta^3\mathcal{D}^3}
		\end{align}
		for any $\eta\leq\frac{c_{\max}}{L_1}$ and any $T\geq\frac{\mathcal{T}}{c_{\max}}$.
	\end{lemma}
   This lemma states that for two points $\u_0$ and $\w_0$, where $\w_0$ lies in the direction of the minimum eigenvector, one of both sequences $\u_t$ and $\w_t$ will result in a further decrement of the objective value and hence escapes the saddle point. To prove the statement of Lemma \ref{funvtion_bound}, we need another two results stated next. In Lemma \ref{bound}, we show that if the function value does not decrease after adding perturbation, then all the iterates belongs to a small ball around $\u_0$. The next Lemma \ref{bound2} shows that if the iterates starting from $\u_0$ get stuck, then the sequence $\w_t$ will result in decreasing the function value and hence escapes the saddle point. Here, $\w_t$ is the sequence obtained after applying SCA algorithm steps to $\w_0$ which by taking a step in the direction of minimum eigenvalue of the Hessian denoted by $\z_1$ and starting from $\u_0$. 
	
	\begin{lemma}\label{bound}
		For any constant $\hat{c}\geq \frac{6}{\left(12+\frac{2\mathcal{D}}{\gamma L}\log\frac{d\kappa}{\delta}\right)^2}$, there exists a constant $c_{\max}$ for any $\delta\in(0,\frac{d\kappa}{e}]$, $f(\cdot)$ satisfying the assumptions, and $\norm{f(\tilde{\x})}\leq G$, for an initial point $\u_0=\tilde{\x}+\zeta$ with $\norm{\u_0-\tilde{\x}}\leq 2r$, let
		\begin{align}
		T:=\min\{\inf_{t}\{t|\tilde{f}_{\u_0}(\u_t)-f(\u_0)\leq-3F\},\hat{c}F\}.
		\end{align} 
		Then for any $\eta\leq c_{\max}/L_1$, it holds that for all $t<T$ we have $\norm{\u_t-\tilde{\x}}\leq Z(L\hat{c})$, where $Z:=48+\frac{8\mathcal{D}}{\gamma L}\log\left(\frac{d\kappa}{\delta}\right)$ with $\mathcal{D}=L_0\left(1+\frac{1}{C}\right)$.
	\end{lemma}
	
	\begin{lemma}\label{bound2}
		There exists $\hat{c}$, $c_{\max}$ for any $\delta\in(0,\frac{d\kappa}{e}]$, $f(\cdot)$ satisfying the assumptions, and $\norm{f(\tilde{\x})}\leq G$, for an initial point $\u_0=\tilde{\x}+\zeta$ and sequence $\{\u_t\},\{\w_t\}$ where $\w_0=\u_0+\mu r \z_1$ with $\mu\in[\frac{\delta}{2\sqrt{2}},1]$, let us define 
		\begin{align}
		T=\min\{\inf_{t}\{t|\tilde{f}_{\w_0}(\w_t)-f(\w_0)\leq-3F\},\hat{c}F\}
		\end{align} 
		then for any $\eta\leq c_{\max}/L_1$, if $\norm{\u_t-\tilde{\x}}\leq Z(L\hat{c})$  for all $t<T$, it holds that $T<\hat{c}F$. 
	\end{lemma}
	
	The proof of Lemmas \ref{function_value} - Lemma \ref{bound2} are provided in the supplementary material available at \cite{supp_paper}. After stating all the intermediate results, we are in position to discuss the proof of Theorem \ref{thm1}. The detailed proof is provided next.\\ 
	
	\begin{myproof}[of Theorem \ref{thm1}]
		Let $c< \min\{c_{\max},2C\}$, the goal of this proof is to achieve a point $\x$ for which it holds that
		\begin{align}\label{condition}
		\!\!\!\!\!\!\norm{\nabla U(\x)}\leq g_\text{th}=\frac{\epsilon\sqrt{c}}{\chi^2}, \ \ \  \ \lambda_{\min}(\nabla^2U(\x))\geq -\sqrt{\epsilon L_2}. 
		\end{align}
		Note that $c\leq 1$, $\chi\geq 1$ as defined in Algorithm \ref{algo3}, which implies that $\frac{\sqrt{c}}{\chi^2}\leq 1$ and hence any $\x$ satisfying \eqref{condition} is an $\epsilon$ second order stationary point. To proceed with the proof, let us start from $\x_0$, if $\x_0$ does not satisfy \eqref{condition}, then there are two possible outcomes
		\begin{itemize}
			\item[1.] $\norm{\nabla U(\x_0)}>g_\text{th}$. This means that the first order stationary point is not reached yet and hence Algorithm 2 will not perform the perturbation step. From Lemma \ref{lem:descent}, we have 
			\begin{align}
			U(\x_1)-U(\x_0)\leq -\frac{\eta}{2}g_\text{th}^2=-\frac{c^2\epsilon^2}{2\chi^4L_1}.
			\end{align}
			\item[2.] $\norm{\nabla U(\x_0)}\leq g_\text{th}$. Under this  condition, a perturbation step is performed by Algorithm 2, then SCA steps are performed for the next $t_\text{th}$ iterations, and then termination condition is checked. If termination condition is not satisfied, it holds that
			\begin{align}\label{bound1}
			U(\x_\text{th})-U(\x_0)\leq -f_\text{th}+\frac{16L_2c^3\mathcal{D}^3}{L_1^3}.\end{align}
		\end{itemize}
		where we have substituted step size $\eta=\frac{c}{L_1}$. Now we select {$c^3=\frac{f_\text{th}}{16L_2\mathcal{D}^3}$ with $s\in(0,1)$ and introduce $s$ as an adjustment parameter to make $c< \min\{c_{\max},2C\}$. Note that we cannot make $s$ arbitrarily small because that would result in s smaller step size $\eta$ which is proportional to $c$. Next, we obtain} \begin{align}
		U(\x_\text{th})-U(\x_0)\leq -(1-s)f_\text{th}.\end{align}
		\
		The above results implies that on an average with each step of the algorithm, the function value reduces by 
		\begin{align}
		\frac{U(\x_\text{th})-U(\x_0)}{t_\text{th}}\leq -\frac{c^3\epsilon^2}{\chi^4L_1}.
		\end{align}
		The function values decreases at least by $\frac{c^3\epsilon^2}{\chi^4L_1}$ on an average after running for $t_\text{th}$ iterations. Since the maximum amount by which the function value may decrease is upper bound by $U(\x_0)-U^\star$, the algorithm must terminate in the following number of iterations
		\begin{align}
		\frac{U(\x_0)-U^\star}{\frac{c^3\epsilon^2}{\chi^4L_1}}=&\frac{\chi^4}{c^3}\cdot\frac{L_1(U(\x_0)-U^\star)}{\epsilon^2}\nonumber 
		\\
		=&\mathcal{O}\left(\frac{L_1(U(\x_0)-U^\star)}{\epsilon^2}\log\left(\frac{dL_1\triangle_U}{\epsilon^2\delta}\right)\right).
		\end{align}
		From Algorithm \ref{algo3}, note that the perturbation is added at iteration $t$ if the gradient is small or $\nabla U(\tilde{x}_t)\leq g_\text{th}$. From Lemma \ref{function_value}, the probability of happening this at each time is at least $1-\delta$ where $\delta=\frac{dL_1}{\sqrt{L_2\epsilon}}\exp(-\chi)$. In other words, the number of times perturbation is added for one run of the algorithm is given by 
		\begin{align}
		\!\!\!\!\!\frac{1}{t_\text{th}}\cdot\frac{\chi^4}{c^3}\cdot\frac{L_1(U(\x_0)-U^\star)}{\epsilon^2}=\frac{\chi^3}{c}\cdot\frac{L_1(U(\x_0)-U^\star)}{\epsilon^2}.
		\end{align}
		Using the union bound, it holds that Lemma \ref{function_value} holds for each of the iterations after adding perturbations. This result makes sure that the Algorithm \ref{algo3} converges to the $\epsilon$ second order stationary point with probability
		\begin{align}
		1-\frac{dL_1}{\sqrt{L_2\epsilon}}&\exp(-\chi)\cdot\frac{\chi^3}{c}\frac{L_1(U(\x_0)-U^\star)}{\epsilon^2}\nonumber \\
		&=1-\frac{\chi^3\exp(-\chi)}{c}\cdot\frac{dL_1(U(\x_0)-U^\star)}{\epsilon^2}.
		\end{align}
		Note that we choose $\mathcal{X}=3\max\{\log\left(\frac{dL_1\triangle_U}{c\epsilon^2}\right),4\}$, which implies that $\chi\geq 12$, and $\chi^3\exp(-\chi)\leq\exp^{-\chi/3}$. Therefore, we have
		\begin{align}
		\frac{\chi^3\exp(-\chi)}{c}\cdot&\frac{dL_1(U(\x_0)-U^\star)}{\epsilon^2}\nonumber \\
		\leq& \frac{\exp^{-\chi/3}}{c}\cdot\frac{dL_1(U(\x_0)-U^\star)}{\epsilon^2}\leq\delta.
		\end{align}
		which completes the proof.
	\end{myproof}

	\section{Conclusion and future directions}
	
    This paper considers an algorithm for solving non-convex optimization formed by perturbation of successive convex approximation  based method. The proposed perturbed successive convex approximation (P-SCA) algorithm is shown to converge to second order stationary point with a rate similar to the vanilla gradient descent (convergence to first order stationary point) up to a constant factor. Hence, the proposed algorithm can escape the saddle point efficiently utilizing the perturbed variant of the successive convex optimization algorithm. To perform the analysis, the proposed algorithm is proved to be a special case of the  standard gradient descent algorithm with the gradient in error, which is called inexact gradient descent in literature. The idea of adding perturbation to the algorithms updates is utilized to escape the saddle points. As a future work, we are interested in considering the constrained version of the proposed algorithm which also escapes the saddle points.

\small 
	\bibliographystyle{named}
	\bibliography{ijcai19}
	
	\setcounter{section}{0}
\normalsize
	\newpage\onecolumn
	\section*{Supplementary Material for ``Escaping Saddle Points with the Successive Convex Approximation Algorithm"}
	
	\titleformat{\section}{\normalfont \Large \bfseries}
	{Appendix\ \thesection: }{2.3ex plus .2ex}{}
	\titlespacing{\subsection}{2em}{*1}{*1}
	
	\section{Proof of   Lemma \ref{function_value}}
	From Lipschitz continuous gradient property, it holds that
	\begin{align}
	U(\x_0)\leq U(\tilde{\x})+\nabla U(\tilde{\x})^T(\x_0-\tilde{\x})+\frac{L_1}{2}\norm{\x_0-\tilde{\x}}^2.
	\end{align}
	From the statement of Lemma \ref{function_value} and Cauchy-Schwarz inequality, it holds that
	\begin{align}
	U(\x_0)-U(\tilde{\x})\leq& \norm{\nabla U(\tilde{\x})}\norm{\bbxi}+\frac{L_1}{2}\norm{\bbxi}^2\\
	\leq & Gr+\frac{L_1}{2}r^2\nonumber\\
	\leq & \frac{GL}{\kappa\log\left(\frac{d\kappa}{\delta}\right)}+\frac{L_1}{2}\left(\frac{L}{\kappa\log\left(\frac{d\kappa}{\delta}\right)}\right)^2 \leq\frac{3}{2}F.\nonumber
	\end{align}
	The above mentioned bound represents the maximum value by which the function value can increase in the worst case after adding the perturbation. Next, following exactly the similar steps as performed in the proof of Lemma 14 in \cite{scutari2017parallel}, we obtain
	\begin{align}
	U(\x_T)-U(\tilde{x})=& U(\x_T)-U(\x_0)+U(\x_0)-U(\tilde{x})\nonumber\\
	\leq & -2.5F+16L_2\eta^3\mathcal{D}^3+1.5F\nonumber \\ \leq& -F+16L_2\eta^3\mathcal{D}^3.
	\end{align}
	
	\section{Proof of  Lemma \ref{funvtion_bound}}
	%
	We consider $\tilde{\x}=0$ without loss of generality, $T^\star=\hat{c}\mathcal{T}$, and $T'$ defined as
	\begin{align}
	T'=\inf_{t}\{t|\tilde{U}_{\u_0}(\u_t)-U(\u_0)\leq -3F\}.
	\end{align}
	In order to prove the lemma statement, let us consider the two cases. 
	
	\textbf{Case I ($T'\leq T^\star$):} From the statement of Lemma \ref{bound}, we have $\norm{\u_{T'-1}}\leq \mathcal{O}(L)$ which implies that
	\begin{align}
	\norm{\u_{T'}}=&\norm{\u_{T'-1}}+\eta\norm{\nabla U(\u_{T'-1})} +\eta\norm{\e_{T'-1}}
	\\
	\leq&\norm{\u_{T'-1}}+\eta \norm{\nabla U(\tilde{\x})}+\eta L_1\norm{\u_{T'-1}} +\eta\norm{\e_{T'-1}}\nonumber
	\\
	\leq&(1+\eta L_1)\hat{c} ZL+ \eta G+\eta\mathcal{D}.\nonumber
	\end{align}
	From the equality $\frac{G\log(\frac{d\kappa}{\delta})}{\gamma}=L$, we can write 
	\begin{align}
	\norm{\u_{T'}}\leq& Z_1L+\eta\mathcal{D}
	\end{align} 
	where $Z_1:=\left((1+\eta L_1)\hat{c}Z+ \eta \frac{\gamma}{\log(\frac{d\kappa}{\delta})}\right)$. By choosing a small $c_{\max}$ and $\eta L_1\leq c_{\max}$, for the Hessian Lipschitz function, we have
	
	\begin{align}
	U(\u_{T'})-U(\u_0)\leq & \nabla U(\u_0)^T(\u_{T'}-\u_0)\\&+\frac{1}{2}(\u_{T'}-\u_0)^T\nabla^2U(\u_0)(\u_{T'}-\u_0)\nonumber\\&+\frac{L_2}{6}\norm{\u_{T'}-\u_0}^3\nonumber
	\end{align}
	Consider the following quadratic approximation 
	\begin{align}
	\tilde{U}_{\mathbf{y}}(\x):=& U(\mathbf{y})+\nabla U(\mathbf{y})^T(\x-\mathbf{y})
	\\&+\frac{1}{2}(\x-\mathbf{y})^T\nabla^2U(\tilde{\x})(\x-\mathbf{y}).\nonumber
	\end{align}
	For a Hessian Lipschitz function, it holds that 
	\begin{align}
	U&(\u_{T'})-U(\u_0)
	\\
	\leq& \nabla U(\u_0)^T(\u_{T'}-\u_0)+\frac{1}{2}(\u_{T'}-\u_0)^T\nabla^2U(\u_0)(\u_{T'}-\u_0)\nonumber 
	\\ 
	&+\frac{L_2}{6}\norm{\u_{T'}-\u_0}^3\nonumber
	\\
	\leq& \nabla U(\u_0)^T(\u_{T'}-\u_0)+\frac{1}{2}(\u_{T'}-\u_0)^T\nabla^2U(\tilde{\x})(\u_{T'}-\u_0)\nonumber\\
	&+\frac{1}{2}(\u_{T'}-\u_0)^T(\nabla^2U(\u_0)-\nabla^2U(\tilde{\x}))(\u_{T'}-\u_0)\nonumber 
	\\ 
	&+\frac{L_2}{6}\norm{\u_{T'}-\u_0}^3\nonumber
	\\
	\leq &\tilde{U}_{\u_0}(\u_{T'})-U(\u_0)+\frac{L_2}{2}\norm{\u_{T'}-\tilde{\x}}\norm{\u_{T'}-\u_0}^2\\
	&+\frac{1}{2}(\u_{T'}-\u_0)^T\nabla^2U(\u_0)(\u_{T'}-\u_0)+\frac{L_2}{6}\norm{\u_{T'}-\u_0}^3\nonumber
	\\ \leq & -3F+\mathcal{O}(L_2L^3)+16L_2\eta^3\mathcal{D}^3\nonumber 
	\\
	=& -3F+\mathcal{O}{(\sqrt{\eta L_1} F)}+16L_2\eta^3\mathcal{D}^3\nonumber
	\\
	\leq & -2.5F {+16L_2\eta^3\mathcal{D}^3}
	\end{align}
	by  selecting $c_{\max}\leq \min\{1,\frac{1}{\hat{c}}\}$ and $\eta<\frac{1}{L_1}$. From Lemma 1, it holds that for any $T\geq\frac{1}{c_{\max}}\mathcal{T}\geq \hat{c}\mathcal{T}=T^\star\geq T'$, we can write
	\begin{align}
	U(\u_{T})-U(\u_0)\leq& U(\u_{T^\star})-U(\u_0)\\ 
	\leq & U(\u_{T'})-U(\u_0) \leq -2.5F+16L_2\eta^3\mathcal{D}^3.\nonumber
	\end{align} 
	
	\textbf{Case II ($T'> T^\star$):} For this case also, we know that $\norm{\u_t}\leq \mathcal{O}(L)+\eta\mathcal{D}$ for all $t\leq T^\star$. Let us define $T''$ as
	\begin{align}
	T''=\inf_{t}\{t\ |\ \tilde{U}_{\w_0}(\w_t)-U(\w_0)\leq -3F \}.
	\end{align}
	From Lemma \ref{bound2}, we have $T''\leq T^\star$. Using the similar arguments as in case 1, we conclude that for all $T\geq \frac{1}{c_{\max}}\mathcal{T}$, it holds that
	\begin{align}
	U(\w_T)-U(\w_0)\leq U(\w_{T^\star})-U(\w_0)\leq -2.5F +16L_2\eta^3\mathcal{D}^3.
	\end{align}
	This concludes the proof.
	

	\section{Proof of Lemma \ref{bound}}
	\begin{myproof}
		The idea here is to show that if function value does not decrease after performing SCA iterations then it must be bounded in a small ball. This  is performed by analyzing the update dynamics of the proposed algorithm via decomposing $d$ dimensional update into two components: one along the subspace $\mathcal{S}$ which is span of significantly negative eigenvalues of the Hessian and second to the subspace compliment to it.
		
		Consider the second order Taylor series approximation of the objective function $U(\x)$ evaluated at $\tilde{\x}$ given by 
		\begin{align}\label{eq:approximation}
		\!\!\!\tilde{U}_{\tilde{\x}}(\x)\!:=\!U(\tilde{\x})\!\!+\!\!\nabla U(\tilde{\x})^T(\x\!\!-\!\!\tilde{\x})\!\!+\!\!\frac{1}{2}(\x-\tilde{\x})^T\mathcal{H}(\x-\tilde{\x})
		\end{align} 
		where $\mathcal{H}=\nabla^2 U(\tilde{\x})$. It is emphasized that the gradient of $U$ can be approximated by gradient of $\tilde{U}_{\tilde{\x}}$ provided $\x$ and $\tilde{\x}$ are close. This result is stated as \cite{nesterov2013introductory}
		\begin{align}
		\norm{\nabla U(\x)-\nabla \tilde{U}_{\tilde{\x}}(\x)}\leq \frac{L_2}{2}\norm{\x-\tilde{\x}}^2
		\end{align} for a $L_2$-Hessian Lipschitz function $U$. Set $\u_0=0$, from the SCA steps in \eqref{SCA-GD}, we can write
		\begin{align}
		\u_{t+1}=&\u_t-\eta\nabla U(\u_t) +\eta\e_t\\
		=&\u_t-\eta\nabla U(0)-\eta\left[\int_{0}^{1}\nabla^2U(\theta\u_t)d\theta\right]\u_t +\eta\e_t\\
		=&\u_t-\eta\nabla U(0)-\eta(\mathcal{H}+\triangle_t)\u_t+\eta \e_t\\
		=&(\mathbf{I}-\eta\mathcal{H}-\eta\triangle_t)\u_t-\eta(\nabla U(0)-\e_t)
		\end{align}
		where $\triangle_t:=\int_{0}^{1}\nabla^2U(\theta\u_t)d\theta-\mathcal{H}$. Next, from the Hessian Lipschitz in \eqref{Hessia_lip}, we have $\norm{\triangle_t}\leq L_{2}(\norm{\u_t}+\norm{\tilde{\x}})$. From the smoothness of the gradient in Assumption 1, we have 
		\begin{align}
		\norm{\nabla U(0)-\e_t}&\leq \norm{\nabla U(0)-\nabla U(\tilde{\x})+\nabla U(\tilde{\x})+\e_t}\\
		&\leq \norm{\nabla U(\tilde{\x})}+L_1\norm{\tilde{\x}}+\mathcal{D}\\
		&\leq G+L_1\cdot 2\frac{L}{\kappa\cdot\log\left(\frac{d\kappa}{\delta}\right)}+\mathcal{D}\nonumber\\
		&\leq 3G+\mathcal{D}.
		\end{align} 
		Next, we will calculate projections of $\u_t$ on to the subspaces $S$ and $S_c$, where $S$ is the subspace of eigen vectors whose corresponding eigen values are less than $-\frac{\gamma}{\hat{c}\log\frac{d\kappa}{\delta}}$. Further, $S_c$ defines the subspace of the remaining eigen vectors. 
		
		Next, from the definition of $T$, we know that
		$\tilde{U}_{\u_0}(\u_t)-U(\u_0)> -3F$. Let $\u_0=0$, we get
		\begin{align}
		-3F<\tilde{U}_{0}(\u_t)-U(0)=\nabla U(0)^T\u_t+\frac{1}{2}\u_t^T\mathcal{H}\u_t 
		\end{align}
		Now decomposing the vector $\u_t$ into $\alpha_t$ (projection on to $S$) and $\beta_t$ (projection on to $S_c$), the updates become
		\begin{align}
		\alpha_{t+1}=&(\mathbf{I}-\eta\mathcal{H})\alpha_t-\eta\mathcal{P}_{S}\triangle_t\u_t-\eta\mathcal{P}_{S}(\nabla U(0)-\e_t)\\
		\beta_{t+1}=&(\mathbf{I}-\eta\mathcal{H})\beta_t-\eta\mathcal{P}_{S^c}\triangle_t\u_t-\eta\mathcal{P}_{S^c}(\nabla U(0)-\e_t)\label{beta}.
		\end{align} Utilizing these definitions, we get
		\begin{align}
		-3F\leq \nabla U(0)^T\u_t-\frac{\gamma}{2}\frac{\norm{\alpha_t}^2}{\hat{c}\log\frac{d\kappa}{\delta}}+\frac{1}{2}\beta_t^T\mathcal{H}\beta_t 
		\end{align}
		where the negative sign comes  for the second term comes from the upper bound on the eigen values in subspace $S$. Note that $\norm{\u_t}^2=\norm{\alpha_t}^2+\norm{\beta_t}^2$, we get
		\begin{align}
		-3F&\leq \nabla U(0)^T\u_t-\frac{\gamma}{2}\frac{\norm{\u_t}^2-\norm{\beta_t}^2}{\hat{c}\log\frac{d\kappa}{\delta}}+\frac{1}{2}\beta_t^T\mathcal{H}\beta_t.
		\end{align}
		After rearranging the terms, we get
		\begin{align}
		\norm{\u_t}^2&\leq\frac{2\hat{c}\log\frac{d\kappa}{\delta}}{\gamma}\left(3F+\nabla U(0)^T\u_t+\frac{1}{2}\beta_t^T\mathcal{H}\beta_t\right)+\norm{\beta_t}^2\nonumber\\
		&\hspace{-1cm}\leq 4\max\!\Big\{\!\frac{6F\hat{c}\log\frac{d\kappa}{\delta}}{\gamma},\frac{6G\hat{c}\log\frac{d\kappa}{\delta}}{\gamma}\norm{\!\u_t\!}\!,\!\frac{\beta_t^T\mathcal{H}\beta_t\hat{c}\log\frac{d\kappa}{\delta}}{\gamma},\norm{\beta_t}^2\!\!\Big\}\nonumber.
		\end{align}
		Taking square root on both sides results in
		\begin{align}\label{eq:main}
		&\norm{\u_t}\\&\!\leq \!4\max\!\Big\{\!\!\sqrt{\frac{6F\hat{c}\log\frac{d\kappa}{\delta}}{\gamma}},\frac{6G\hat{c}\log\frac{d\kappa}{\delta}}{\gamma},\sqrt{\frac{\beta_t^T\mathcal{H}\beta_t\hat{c}\log\frac{d\kappa}{\delta}}{\gamma}},\norm{\beta_t}\!\!\Big\}.\nonumber
		\end{align}
		In order to prove the lemma statement 
		\begin{align}\label{lemma_statement}
		\norm{\u_t}\leq \hat{c}LZ
		\end{align}
		we proceed by induction. It holds trivially for $t=0$ since $\u_0=0$. Assume that \eqref{lemma_statement} holds for $\tau\leq t$, and then it remains to show that \eqref{lemma_statement} holds for $t+1<T$. It is sufficient to bound the last two terms of \eqref{eq:main} which are $\beta_t^T\mathcal{H}\beta_t$ and $\norm{\beta_t}$. Consider the update equation for $\beta_t$ from \eqref{beta}, we have 
		\begin{align}\label{eq:beta}
		\beta_t=&(\mathbf{I}-\eta\mathcal{H})\beta_t+\eta\delta_t
		\end{align}
		where $\norm{\delta_t}$ is bounded as
		\begin{align}
		\norm{\delta_t}&\leq\norm{\triangle_t}\norm{\u_t}+\norm{\nabla U(0)-\e_t}\\
		&\leq  L_{2}(\norm{\u_t}+\norm{\tilde{\x}})\norm{\u_t}+\norm{\nabla U(0)-\e_t}\\
		&\leq L_{2}\hat{c}Z\left(\hat{c}Z\!+\frac{2}{\kappa\cdot \log(\frac{d\kappa}{\delta})}\right)L^2\!\!+\!3G\!+\mathcal{D}\nonumber
		\\
		&\leq\left[\hat{c}Z(\hat{c}Z+2)\sqrt{\eta{L_1}}+3\right]G+\mathcal{D} \\
		&\leq 4G+\mathcal{D}
		\end{align}
		where the last step follows by choosing $c_{\max}\leq\frac{1}{Z\hat{c}(Z\hat{c}+2)}$ and step size {$\eta<c_{\max}/L$}.
		
		\textbf{Bounding $\norm{\beta_{t+1}}$:} We have
		\begin{align}
		\norm{\beta_{t+1}}\leq \left(1+\frac{\eta\gamma}{\hat{c}\log\left(\frac{d\kappa}{\delta}\right)}\right)\norm{\beta_t}+\eta(4G+\mathcal{D})
		\end{align}
		Applying recursively, we get
		\begin{align}
		\norm{\beta_{t+1}}&\leq \sum\limits_{\tau=0}^{t}(4G+\mathcal{D})\left(1+\frac{\eta\gamma}{\hat{c}\log\left(\frac{d\kappa}{\delta}\right)}\right)^\tau \eta \nonumber
		\\
		&\leq 3(4G+\mathcal{D})\eta T  
		\\
		&\leq (12 L+\eta\mathcal{D}\mathcal{T})\hat{c}.\label{beta_bound_1}
		\end{align}
		Now, from the definition we have $T\leq\hat{c}F$, so that $\left(1+\frac{\eta\gamma}{\hat{c}\log\left(\frac{d\kappa}{\delta}\right)}\right)^T\leq 3$. 
		
		\textbf{Bounding $\beta_{t+1}^T\mathcal{H}\beta_{t+1}:$} From \eqref{eq:beta}, we can write
		\begin{align}
		\beta_t=\sum\limits_{\tau=0}^{t-1}(\mathbf{I}-\eta\mathcal{H})^\tau\delta_\tau,
		\end{align}
		which implies that
		\begin{align}
		\!\beta_{t+1}^T\mathcal{H}\beta_{t+1}&\!=\!\eta^2\!\!\sum\limits_{\tau_1=0}^{t}\sum\limits_{\tau_2=0}^{t}\delta_{\tau_1}^T(\mathbf{I}\!-\!\eta\mathcal{H})^{\tau_1}\mathcal{H}(\mathbf{I}-\eta\mathcal{H})^{\tau_2}\delta_{\tau_2}\\
		&\hspace{-18mm}\leq \eta^2\sum\limits_{\tau_1=0}^{t}\sum\limits_{\tau_2=0}^{t}\norm{\delta_{\tau_1}}\norm{(\mathbf{I}-\eta\mathcal{H})^{\tau_1}\mathcal{H}(\mathbf{I}-\eta\mathcal{H})^{\tau_2}}\norm{\delta_{\tau_2}}\\
		&\hspace{-18mm}\leq (4G+\mathcal{D})^2\eta^2\sum\limits_{\tau_1=0}^{t}\sum\limits_{\tau_2=0}^{t}\norm{(\mathbf{I}-\eta\mathcal{H})^{\tau_1}\mathcal{H}(\mathbf{I}-\eta\mathcal{H})^{\tau_2}}.
		\end{align}
		If $\{\lambda_i\}$ denotes eigen values of $\mathcal{H}$, then $\lambda_i(1-\eta\lambda_i^{\tau_1+\tau_2})$ denotes the corresponding eigen value of $(\mathbf{I}-\eta\mathcal{H})^{\tau_1}\mathcal{H}(\mathbf{I}-\eta\mathcal{H})^{\tau_2}$. The maxima is achieved for $\lambda_{i,t}^\star=\frac{1}{(1+t)\eta}$. Note that the function $\lambda_i(1-\eta\lambda_i^{\tau_1+\tau_2})$ is monotonically increasing between $(-\infty,\lambda_{i,t}^\star]$. Hence,
		\begin{align}
		\norm{(\mathbf{I}-\eta\mathcal{H})^{\tau_1}\mathcal{H}(\mathbf{I}-\eta\mathcal{H})^{\tau_2}}&=\max_i\lambda_i(1-\eta\lambda_i)^{\tau_1+\tau_2}\nonumber\\
		&\leq\hat{\lambda}(1-\eta\hat{\lambda})^{\tau_1+\tau_2}\nonumber\\
		&\leq \frac{1}{(1+\tau_1+\tau_2)\eta}
		\end{align}
		where $\hat{\lambda}=\min\{L_1,\lambda_{\tau_1+\tau_2}^\star\}$. Hence, we get
		\begin{align}
		\beta_{t+1}^T\mathcal{H}\beta_{t+1}&\leq 2(4G+\mathcal{D})^2\eta T\nonumber
		\\
		&\leq 4(4G^2\eta \hat{c}\mathcal{T}+\eta \hat{c}\mathcal{T}\mathcal{D}^2)\nonumber
		\\&\leq 2(4\hat{c}L^2\gamma \log^{-1}\left(\frac{d\kappa}{\delta}\right)+\eta \hat{c}\mathcal{T}\mathcal{D}^2)\label{beta_bound_2}
		\end{align}
		which follows from 
		\begin{align}
		\sum\limits_{\tau_1=0}^{t}\sum\limits_{\tau_2=0}^{t}\frac{1}{(1+\tau_1+\tau_2)\eta}&=\sum\limits_{\tau=0}^{2t}\min\{1+\tau,2t+1-\tau\}\cdot\frac{1}{1+\tau}\nonumber\\&\leq 2t+1<2T.
		\end{align}
		Finally, substituting the upper bounds in \eqref{beta_bound_1} and \eqref{beta_bound_2} into the right hand side of \eqref{eq:main}, and further simplifying the bounds, we get
		\begin{align}\label{eq:main2}
		\norm{\u_t}&\leq ZL\hat{c}.
		\end{align}
		where $Z:=48+\frac{8\mathcal{D}}{\gamma L}\log\left(\frac{d\kappa}{\delta}\right)$ is a constant and $\hat{c}$ is selected such that $\hat{c}\geq \frac{6}{\left(12+\frac{2\mathcal{D}}{\gamma L}\log\frac{d\kappa}{\delta}\right)^2}$. This concludes the proof.
	\end{myproof}
	
	\section{Proof of Lemma \ref{bound2}}
	
	\begin{myproof}
		\normalfont
		Consider the update equation for $\w_t$, we have
		\begin{align}
		\u_{t+1}+\bbv_{t+1}=&\w_{t+1}-\eta\nabla U(\w_t)+\eta\e_t\\
		&\hspace{-2cm}=\u_{t}+\bbv_{t}-\eta\nabla U(\u_{t}+\v_{t})\\
		&\hspace{-2cm}=\u_{t}+\bbv_{t}-\eta\nabla U(\u_{t})-\eta\left[\int_{0}^{1}\nabla^2U(\u_t+\theta\v_t)d\theta\right]\v_t+\eta\e_t\nonumber\\
		&\hspace{-2cm}=\u_{t}+\bbv_{t}-\eta\nabla U(\u_{t})-\eta(\mathcal{H}+\triangle_t')\v_t+\eta\e_t\\
		&\hspace{-2cm}=\u_{t}-\eta\nabla U(\u_{t})+\eta\e_t+\eta(\mathbf{I}-\eta\mathcal{H}-\triangle_t')\v_t
		\end{align}
		From now onwards, the proof is exactly similar to that of Lemma 17 in \cite{scutari2017parallel} but provided here for completeness. As in the earlier proofs, , we have $\norm{\triangle_t'}\leq L_{2}(\norm{\u_t}+\norm{\v_t}+\norm{\tilde{\x}})$. Now, we can write the update for $\v_{t+1}$ as 
		\begin{align}\label{v_t_udpate}
		\v_{t+1}=(\mathbf{I}-\eta\mathcal{H}-\triangle_t')\v_t
		\end{align}
		Note that we have $\norm{\w_{0}-\tilde{\x}}\leq\norm{\u_0-\tilde{\x}}+\norm{\v_{0}}\leq\frac{L}{\kappa\cdot\log\left(\frac{d\kappa}{\delta}\right)}$, from Lemma \eqref{bound}, we have $\norm{\w_t}\leq ZL\hat{c}$ for all $t\leq T$. From the condition in Lemma \ref{bound2}, we note that $\norm{\u_t}\leq ZL\hat{c}$ for all $t<T$. This implies that \begin{align}\label{bound_last_0}
		\norm{\v_t}\leq\norm{\u_t}+\norm{\w_t}\leq 2ZL\hat{c}
		\end{align} for all $t<T$. From here, we can derive that 
		\begin{align}
		\norm{\triangle_t'}\leq& L_{2}(\norm{\u_t}+\norm{\v_t}+\norm{\tilde{\x}})\\
		\leq &L_2\left(2ZL\hat{c}+\frac{L}{\kappa\log\frac{d\kappa}{\delta}}\right)\\
		\leq & L_2L\left(2Z\hat{c}+1\right).
		\end{align} 
		
		Next, we develop a lower bound on $\norm{\v_t}$. Let $\psi_t$ be the norm of the projection of $\v_t$ onto {$\z_1$ direction} and  $\phi_t$ be the norm of the projection of $\v_t$ on to the remaining subspace. From \eqref{v_t_udpate}, we can write
		\begin{align}
		\psi_t\geq&(1+\gamma\eta)\psi_t-\mu\sqrt{\psi_t^2+\phi_t^2}\\
		\phi_t \leq&(1+\gamma\eta)\Psi_t+\mu\sqrt{\psi_t^2+\phi_t^2}
		\end{align}
		where $\mu=\eta L_2L\left(2ZL\hat{c}+1\right)$. By induction, next we prove that the following inequality holds for all $t<T$
		\begin{align}\label{induction}
		\phi_t\leq 4\mu t\psi_t.
		\end{align} 
		For $t=1$ the base case of the induction holds for $t=0$ from the definition of $\v_0$. Now, let us assume that \eqref{induction} holds for $t\leq T$,  we need to show that it holds for $t+1\leq T$, we have
		\begin{align}
		4\mu(t+1)\psi_t\geq& 4\mu(t+1)\left((1+\gamma\eta)\psi_t-\mu\sqrt{\psi_t^2+\phi_t^2}\right)\\
		\phi_{t+1}\leq &4\mu t(1+\gamma\eta)\psi_t+\mu\sqrt{\psi_t^2+\phi_t^2}.
		\end{align}
		Next, we only need to show that
		\begin{align}
		(1+4\mu(t+1))\sqrt{\psi_t^2+\phi_t^2}\leq 4(1+\gamma\eta)\psi_t.
		\end{align}
		By selecting $\sqrt{c_{\max}}\leq\frac{1}{2Z\hat{c}+1}\min\{\frac{1}{2\sqrt{2}},\frac{1}{4\hat{c}}\}$ and $\eta\leq\frac{c_{\max}}{L_1}$, we get
		\begin{align}
		4\mu(t+1)\leq 4\mu T\leq & 4\eta L_2 L(2Z\hat{c}+1)\hat{c}\mathcal{T}\nonumber\\=&4\sqrt{\eta L_1}(2Z\hat{c}+1)\hat{c}\leq 1.
		\end{align}
		This implies that
		\begin{align}
		4(1+\gamma\eta)\psi_t\geq 4\psi_t\leq 2\sqrt{2\phi_t^2}\geq (1+4\mu(t+1))\sqrt{\psi_t^2+\phi_t^2}
		\end{align}
		which proves the induction. From \eqref{induction}, we have $\phi_t\leq 4\mu t\psi_t\leq \psi_t$, which gives
		\begin{align}\label{bound_last}
		\psi_{t+1}\geq (1+\gamma\eta)\psi_t-\sqrt{2}\mu\psi_t\geq\left(1+\frac{\gamma\eta}{2}\right)\psi_t
		\end{align}
		where the last inequality follows from $\mu=\eta L_2 L(2ZL\hat{c}+1)\leq \sqrt{c_{\max}}(2ZL\hat{c}+1)\gamma\eta \log^{-1}\frac{d\kappa}{\delta}<\frac{\gamma\eta}{2\sqrt{2}}$.
		
		From \eqref{bound_last_0} and \eqref{bound_last}, we obtain for all $t<T$
		\begin{align}
		2ZL\hat{c}\geq& \norm{\v_t}\geq \psi_t\geq \left(1+\frac{\gamma\eta}{2}\right)^t\psi_0\nonumber\\
		=&\left(1+\frac{\gamma\eta}{2}\right)^tc_0\frac{L}{\kappa}\log^{-1}\frac{d\kappa}{\delta}\nonumber\\
		\geq &\left(1+\frac{\gamma\eta}{2}\right)^t\frac{\delta}{2\sqrt{d}}\frac{L}{\kappa}\log^{-1}\left(\frac{d\kappa}{\delta}\right).
		\end{align}
		The above result implies that
		\begin{align}
		T\leq&\frac{1}{2}\frac{\log\left(4Z\frac{\kappa\sqrt{d}}{\delta}\cdot\hat{c}\log\left(\frac{d\kappa}{\delta}\right)\right)}{\log\left(1+\frac{\gamma\eta}{2}\right)}\nonumber\\
		\leq & \frac{1}{2}\frac{\log\left(4Z\frac{\kappa\sqrt{d}}{\delta}\cdot\hat{c}\log\left(\frac{d\kappa}{\delta}\right)\right)}{\gamma\eta}\nonumber\\
		\leq & (2+\log(4Z\hat{c}))\mathcal{T}.
		\end{align}
		The inequality follows from the selection $\delta\in(0,\frac{d\kappa}{\delta}]$ and we have $\log\left(\frac{d\kappa}{\delta}\right)\geq 1$. We choose $\hat{c}$ large enough such that $2+\log(4Z\hat{c})\leq \hat{c}$, we get $T<\hat{c}\mathcal{T}$, which concludes the proof.
	\end{myproof}
	
	%
	%
	%
	%
	%

\end{document}